\documentclass[12pt]{amsart}

\usepackage{latexsym}
\usepackage[dvips]{graphics}

\usepackage{amsmath,amssymb,amscd,amsthm,amsxtra}

\newtheorem{theorem}{Theorem}[section]
\newtheorem{proposition}[theorem]{Proposition}
\newtheorem{corollary}[theorem]{Corollary}
\newtheorem{lemma}[theorem]{Lemma}
\theoremstyle{definition}
\theoremstyle{remark}

\newcommand{\q}{\quad}
\newcommand{\qq}{\quad\quad}
\newcommand{\qqq}{\quad\quad\quad}
\newcommand{\qqqq}{\quad\quad\quad\quad}

\newcommand{\nf}{\infty}

\newcommand{\al}{\alpha}

\newcommand{\ga}{\gamma}

\newcommand{\de}{\delta}

\newcommand{\ve}{\varepsilon}

\newcommand{\la}{\lambda}

\newcommand{\rn}{{\mathbf R}^n}

\newcommand{\lp}{L^{p}}

\newcommand{\li}{L^{\infty}}

\newcommand{\cd}{\mathcal D}

\newcommand{\dint}{\displaystyle\int}

\newcommand{\f}{\frac}

\newcommand{\p}{\partial}

\newcommand{\lab}{\label}
\newcommand{\beq}{\begin{equation}}
\newcommand{\eeq}{\end{equation}}
\newcommand{\beqn}{\begin{equation*}}
\newcommand{\eeqn}{\end{equation*}}
\newcommand{\bp}{\begin{proof}}
\newcommand{\ep}{\end{proof}}
\newcommand{\bprop}{\begin{proposition}}
\newcommand{\eprop}{\end{proposition}}
\newcommand{\bt}{\begin{theorem}}
\newcommand{\et}{\end{theorem}}
\newcommand{\bc}{\begin{corollary}}
\newcommand{\ec}{\end{corollary}}
\newcommand{\bl}{\begin{lemma}}
\newcommand{\el}{\end{lemma}}

\newcommand{\cz}{Calder\'on-Zygmund }
\newcommand{\mcz}{multilinear Calder\'on-Zygmund }

\begin{document}

\title
[multilinear operators on Hardy spaces]
{multilinear Calder\'on-Zygmund operators on Hardy spaces}

\author{Loukas Grafakos}

\address{
Loukas Grafakos\\
Department of Mathematics\\
University of Missouri\\
Columbia, MO 65211, USA}

\email{loukas@math.missouri.edu}

\author{Nigel Kalton}

\address{
Nigel Kalton\\
Department of Mathematics\\
University of Missouri\\
Columbia, MO 65211, USA}

\email{nigel@math.missouri.edu}

\thanks{The second author was supported by NSF grant DMS-98700027}

\date{\today}



\begin{abstract}
It is shown that multilinear Calder\'on-Zygmund  operators
are bounded on products of Hardy spaces.
\end{abstract}

\maketitle

\section{Introduction}\label{section1}
The study of multilinear singular integral
operators has recently received
increasing attention.  In analogy to the linear theory, the class of
multilinear singular integrals with standard \cz kernels
provide the foundation and   starting point of
investigation  of the theory.
The class of \mcz operators
was introduced and first investigated by Coifman and Meyer \cite{CM1},
\cite{CM2}, \cite{CM3}, and  was later
systematically studied by Grafakos and Torres  \cite{GT}.

In this article we take up the issue of boundedness of \mcz operators
on products of Hardy spaces.  As in the linear theory, a certain amount of extra
smoothness is required for these operators to have such boundedness properties.
We will assume that $K(y_0,y_1,\dots , y_m)$ is a function
defined  {\em away from the  diagonal}
$y_0=y_1=\dots =y_m$  in $(\rn)^{m+1}$ which satisfies the following
estimates
\begin{equation}\lab{assump}
\q \big|\p_{y_0}^{\al_0} \dots \p_{y_m}^{\al_m}
K(y_0,y_1,\dots , y_m) \big|\le
\f{ A_\al }{\big(\sum\limits_{k,l=0}^m |y_k-y_l|\big)^{mn+|\al|}},
\qq \text{for all $|\al|\le N$,}
\end{equation}
where $\al =(\al_0,\dots , \al_m)$ is an ordered set   of $n$-tuples of
nonnegative integers, $|\al|=|\al_0|+\dots +|\al_m|$, where $|\al_j|$ is
the order of each multiindex $\al_j$, and
$N$ is a large integer to  be determined later.  We will call
such functions $K$ multilinear standard kernels. We
assume throughout that $T$ is a weakly continuous
$m$-linear operator defined on
products  of test functions such that
for some multilinear standard kernel  $K$,
the integral representation below is valid
\begin{equation}\lab{defT}
 T(f_1, \dots  , f_m)(x)= \int_{\rn}\dots\int_{\rn}  K(x,y_1,\dots , y_m)
\prod_{j=1}^m f_j(y_j) \, d y_1, \dots dy_m,
\end{equation}
whenever $f_j$ are smooth functions with compact support and
$x \notin \cap_{j=1}^m
\text{supp} f_j$. In the case $m=1$  conditions (\ref{assump}) are called
standard estimates and operators given by (\ref{defT}) are called
\cz if they are bounded from $L^2(\rn)$ to $L^2(\rn)$.
We will adopt the same terminology in the multilinear case and
call $T$ a \mcz operator if it is
associated to a multilinear standard kernel as in
(\ref{defT}) and   has a bounded extension from
a product of  some $L^{q_j}$ spaces into another $L^q$ space with
$1/q=1/q_1+\dots + 1/q_m$. If this is the case, it was shown in \cite{GT}
that these operators map any other product of Lebesgue spaces
$\prod_{j=1}^mL^{p_j}(\rn)$ with $p_j>1$ into the corresponding $\lp$ space.

When $m=1$, bounded extensions  for \cz operators on the Hardy spaces $H^p$ were
obtained   by Fefferman and Stein \cite{FS}.
Here $H^p=H^p(\rn)$ denotes the real Hardy space
in \cite{FS} defined for $0<p\le 1$.
In this note we provide   analogous bounded extensions
for \mcz on products of Hardy spaces. The following theorem is our
main result:

\begin{theorem}\lab{main}
Let   $1<q_1,\dots , q_m,q<\nf$ be fixed indices satisfying
\begin{equation}\lab{ind1}
\f{1}{q_1}+\dots +\f{1}{q_m}=\f{1}{q}
\end{equation}
and let $0<p_1,\dots , p_m,p\le 1$ be   real numbers satisfying
\begin{equation}\lab{ind2}
 \f{1}{p_1}+\dots +\f{1}{p_m}=\f{1}{p}.
\end{equation}
Suppose that $K$ satisfies (\ref{assump}) with $N=[n(1/p-1)]$.
Let $T$ be related to $K$ as in (\ref{defT}) and assume that $T$
admits an extension that
maps $L^{q_1}(\rn)\times \dots \times L^{q_m}(\rn)$ into
$L^q(\rn)$ with norm $B$.  Then
$T$ extends to a bounded operator from
$H^{p_1}(\rn)\times \dots \times H^{p_m}(\rn)$ into
$L^p(\rn)$ which satisfies the norm estimate
$$
\|T\|_{H^{p_1} \times \dots \times H^{p_m}\to L^p} \le C
\big(B+\sum_{|\al |\le N+1} A_\al \big),
$$
for some   constant $C=C(n,p_j,q_j)$.
\end{theorem}

\section{The proof of Theorem \ref{main}}\label{section2}

\begin{proof}
We prove the theorem using the atomic decomposition of $H^p$.  See
Coifman \cite{coif} and Latter \cite{latter}.
Since finite sums of atoms are dense in $H^p$ we
will work with such sums   and we will
obtain estimates independent of the number of terms in each sum.
The general case will follow by a simple density argument.
Write each $f_j$, $1\le j\le m$ as a finite sum of $H^{p_j}$-atoms
$f_j=\sum_{k} \la_{j,k} a_{j,k}$, where $a_{j,k}$ are $H^{p_j}$-atoms.
This means that the $a_{j,k}$'s  are functions  supported in cubes $Q_{j,k}$
and  satisfy  the  properties
\begin{align}
&   |   a_{j,k} |   \le   |Q_{j,k}|^{-\f{1}{p_j} } \lab{12121}  \\
&\int_{Q_{j,k}} x^\ga a_{j,k}(x) \, dx  =0 \lab{gggg},
\end{align}
for all  $|\ga |\le [n(1/p_j-1)]. $  By the theory of
$H^p$ spaces, see \cite{stein-new} page 112, we can take the atoms $ a_{j,k}$
to have vanishing moments up to any  large  fixed specified integer.
In this article we   will assume that all the $a_{j,k}$'s  satisfy
(\ref{gggg}) for all $|\ga | \le [n(1/p -1)]$.
 For a cube $Q$, let $Q^*$ denote   the cube with
the same center and $2\sqrt{n}$ its   side length,
i.e. $l(Q^*)=2\sqrt{n}l(Q)$.  Using multilinearity  we write
\begin{equation}\lab{eee}
T(f_1,\dots , f_m)(x)= \sum_{k_1}\dots \sum_{k_m} \la_{1,k_1}\dots \la_{m,k_m}
T(a_{1,k_1},\dots , a_{m,k_m})(x).
\end{equation}
We now fix  $k_1,\dots ,k_m$ and $x\in \rn$ and we consider the following
cases:
\newline\noindent Case $1$: $x\in Q_{1,k_1}^*\cap \dots \cap
Q_{m,k_m}^*$
\newline\noindent Case  $2$: $x$ lies  in the complement of at least one
of the cubes  $Q_{j,k_j}^*$.

Let us begin with case 2. We fix   $1 \le r \le m$ and without loss of
generality (by permuting the indices) we assume   that
$x\in Q_{r+1,k_{r+1}}^*\cap \dots \cap Q_{m,k_m}^*$ and
that
$x\notin Q_{1,k_1}^*\cup\dots \cup Q_{r,k_r}^*$. Assume without loss of
generality  that the side length of the cube $Q_{1,k_1}$ is the smallest among
the  side lengths of the cubes $Q_{1,k_1},\dots , Q_{r,k_r}$.
Let  $c_{j,k_j}$ be the center of the cube $Q_{j,k_j}$.
Since $a_{1,k_1}$ has zero vanishing moments up to   order $N=[n(1/p_1-1)]$,
we can subtract the Taylor polynomial $  P^{N}_{c_{1,k_1}}(x,\,\cdot\, ,y_2,
\dots , y_m)$ of  the function $  K(x,\,\cdot\, , y_2, \dots
, y_m)$ at  the point $c_{1,k_1}$  to obtain
\begin{align*}
& T(a_{1,k_1}, \dots , a_{m,k_m})(x) \\
=& \int_{(\rn)^{m-1}} \prod_{j=2}^m a_{j,k_j}(y_j) \int_{\rn}a_{1,k_1}(y_1)
\big[K(x,y_1, \dots , y_m )-P^{N}_{c_{1,k_1}}(x,y_1,y_2,
\dots , y_m)\big]\, d\vec y
\\ =&\int_{(\rn)^{m-1}} \prod_{j=2}^m a_{j,k_j}(y_j) \int_{\rn}a_{1,k_1}(y_1)
\!\! \sum_{|\ga |=N+1}\!\! (\p^\ga_{y_1}K)(x,\xi, y_2, \dots , y_m)
\f{(y_1\! -\! c_{1,k_1})^\ga }{\ga !} \, d\vec y,
\end{align*}
for some $\xi $ on the line segment joining $y_1$ to $c_{1,k_1}$ by Taylor's
theorem.
We have $|x-\xi| \ge |x-c_{1,k_1}| -|\xi -c_{1,k_1}|
\ge  |x-c_{1,k_1}| - \f12 \sqrt{n} \, l(Q_{1,k_1}) \ge \f12 |x-c_{1,k_1}|$,
since $x\notin Q_{1,k_1}^*= 2\sqrt{n} Q_{1,k_1}$.  Similarly we
obtain $|x-y_j| \ge \f12 |x-c_{j,k_j}|$ for $j\in \{2,3,\dots , r\}$.
Set
\begin{equation}\lab{defa}
A=   \sum_{|\ga |\le N+1} A_\ga
\end{equation}
 and note that
$A\ge \sum_{|\ga |= N +1} A_\ga$.
The estimates for the kernel $K$ and the size estimates for the
atoms give the
following pointwise bound for the  expression above:
\begin{align*}
  &\underbrace{\int_{\rn}\!\! \dots \! \int_{\rn}}_{m-r\,\,\text{times}}
\bigg(\sum_{|\ga |= N +1} A_\ga \bigg)  \bigg(
\dint_{Q_{1,k_1}} |a_{1,k_1}(y_1) |\,  |y_1-c_{1,k_1}|^{N+1} dy_1 \bigg)  \\
 &\qqqq\f{ |Q_{2,k_{2}} |^{1-\f{1}{p_2}}\dots  |Q_{r,k_{r}} |^{1-\f{1}{p_r}}
|Q_{r+1,k_{r+1}} |^{ - \f{1}{p_{r+1}}}\dots  |Q_{m,k_{m}} |^{  -\f{1}{p_m}}
 dy_{m}\dots dy_{r+1}
 }{\big(\f12 |x-c_{1,k_1}|+  \dots +\f12 |x-c_{r,k_r}| +|x-y_{r+1}|+
\dots + |x-y_m| \big)^{nm+N+1 }} .
\end{align*}
Integrating the above over $y_m \in \rn$, $\dots$, $y_{r+1}\in \rn$ we
obtain that the expression above is bounded by a constant
multiple of
$$
 \f{ A  |Q_{2,k_{2}} |^{1-\f{1}{p_2}}\dots  |Q_{r,k_{r}} |^{1-\f{1}{p_r}}
\dint_{Q_{1,k_1}} |a_{1,k_1}(y_1) |\,  |y_1-c_{1,k_1}|^{N+1} dy_1
 }{\big(\f12 |x-c_{1,k_1}|+  \dots +\f12 |x-c_{r,k_r}|
   \big)^{nm+N+1-n(m-r)}
 |Q_{r+1,k_{r+1}} |^{  \f{1}{p_{r+1}}}\dots  |Q_{m,k_{m}} |^{  \f{1}{p_m}}} .
$$
But the  integral above is easily seen to be controlled by  a
constant multiple of $|Q_{1,k_1}|^{1-\f{1}{p_2}+\f{N+1}{n}}$.
Since the cube $Q_{1,k_1}$ was picked to have
the smallest  size among the $Q_{1,k_1}, \dots , Q_{r, k_r}$, the
expression above    is
bounded by a constant multiple of
\begin{align*}
&A  \prod_{j=1}^r \f{  |Q_{j,k_{j}}
|^{1-\f{1}{p_j} +\f{N+1}{nr}}   }{
\big(|x-c_{j,k_j}|+l(Q_{j,k_{j}}) \big)^{n+\f{  N+1}{r} }  }\,
 |Q_{r+1,k_{r+1}}|^{  -\f{1}{p_{r+1}} }\dots | Q_{m,k_{m}}|^{  -\f{1}{p_m} } \\
\le &C\,A
\prod_{j=1}^m \f{  |Q_{j,k_{j}}
|^{1-\f{1}{p_j} +\f{N+1}{nr}}   }{
\big(|x-c_{j,k_j}|+l(Q_{j,k_{j}}) \big)^{n+\f{ N+1}{r} }  },
\end{align*}
 since
$x\in  Q_{r+1,k_{r+1}}^*\cap \dots \cap Q_{m,k_m}^*   $.

Summing over  all possible $1 \le r \le m$ and all
possible combinations of subsets of
$\{1, \dots , m\}$ of size $r$   we   obtain   the pointwise estimate
\begin{equation}\lab{rrr}
|T(a_{1,k_1}, \dots , a_{m,k_m})(x)| \le C  A
\prod_{j=1}^m \f{  |Q_{j,k_{j}}
|^{1-\f{1}{p_j} +\f{N+1}{nr}}   }{
\big(|x-c_{j,k_j}|+l(Q_{j,k_{j}}) \big)^{n+\f{  N+1}{r} }  }\,
\end{equation}
for  all $x$ which  belong to the complement of at least one
$Q_{j,k_{j}}^*  $ (case 2).

Now using (\ref{eee}) and (\ref{rrr}) we obtain
$$
 | T(f_1, \dots , f_m) (x)|    \le     G_1(x) +G_2(x),
$$
where $G_1(x)$ and $G_2(x)$ correspond to cases $1$ and $2$
respectively and are given by
\begin{align*}
G_1(x) & =
\sum_{k_1}\dots \sum_{k_m} |\la_{1, k_1}| \dots |\la_{m, k_m}|
| T(a_{1,k_1}, \dots , a_{m,k_m}) (x)| \chi_{Q_{1,k_1}^*\cap \dots \cap
Q_{m,k_m}^*}(x)
\\ G_2(x) & =C\, A \prod_{j=1}^m \bigg(\sum_{k_j} |\la_{j, k_j}| \f{
|Q_{j,k_{j}} |^{1-\f{1}{p_j} +\f{N+1}{nr}}   }{
\big(|x-c_{j,k_j}|+l(Q_{j,k_{j}}) \big)^{n+\f{  N+1}{r} }  } \bigg)  .
\end{align*}

Applying H\"older's inequality with exponents $p_1, \dots , p_m$ and $p$ we
obtain the  estimate
\begin{align}\begin{split}\lab{two}
 \|G_2 \|_{\lp} \le
&C A  \prod_{j=1}^m  \bigg\|\sum_{k_j} |\la_{j, k_j}| \f{  |Q_{j,k_{j}}
|^{1-\f{1}{p_j} +\f{N+1}{nr}}   }{
\big(|x-c_{j,k_j}|+l(Q_{j,k_{j}}) \big)^{n+\f{  N+1}{r} }  } \bigg\|_{L^{p_j}}\\
\le & C' A   \prod_{j=1}^m \big(\sum_{k_j} |\la_{j,
k_j}|^{p_j}\big)^{\f{1}{p_j}}
\le C' A  \prod_{j=1}^m  \|f_j\|_{H^{p_j}},
\end{split}\end{align}
where we used the $p_j$-subadditivity of the   $L^{p_j}$ quasi-norm  and
the easy fact that  the  functions
$$
\f{  |Q_{j,k_{j}}
|^{1-\f{1}{p_j} +\f{N+1}{nr}}   }{
\big(|x-c_{j,k_j}|+l(Q_{j,k_{j}}) \big)^{n+\f{  N+1}{r} }  }
$$
have $L^{p_j}$ norms bounded by constants.

We now turn our attention to case 1. Here we will  show that
\begin{equation}\lab{one}
\|G_1 \|_{\lp} \le
C (A+B) \prod_{j=1}^m  \|f_j\|_{H^{p_j}}.
\end{equation}
To prove (\ref{one}) we will need the following lemma whose
proof we postpone until the next section.

\begin{lemma}\lab{L1} Let $0<p\le 1.$
Then  there is a constant $C(p)$ such that  for all finite collections
of cubes $\{Q_k\}_{k=1}^K$ in $\rn$
and all  nonnegative integrable functions
$g_k $ with  $\text{supp } g_k\subset Q_k $  we have
$$
\big\|\sum_{k=1}^K g_k \big\|_{\lp} \le C(p)\,  \bigg\|\sum_{k=1}^K \bigg(
 \f{1}{|Q_k| }\int_{Q_k} g_k(x)\, dx
\bigg) \chi_{Q^*_k}  \bigg\|_{L^p} .
$$
\end{lemma}

We momentarily assume Lemma \ref{L1} and we prove (\ref{one}). Using
the assumption that $T$ maps $L^{q_1}\times \dots \times L^{q_m}$,
it was proved  in \cite{GT}   that $T$  maps
  all possible combinations   of products
$$
\li_c \times \dots  \times \li_c \times  L^2 \times
\li_c \times \dots  \times \li_c
$$
into $L^2$
with norm at most a multiple of $A+B$.  $L^\nf_c$ denotes here the
space of all $\li$ functions with compact support.

Now fix
atoms $a_{1,k_1}$, $\dots$, $a_{m,k_m}$
supported in cubes $Q_{1,k_1}$, $\dots$,    $Q_{m,k_m}$ respectively.
Assume that $Q_{1,k_1}^*\cap\dots \cap    Q_{m,k_m}^* \neq \emptyset$,
otherwise there is nothing to prove.
Since $Q_{1,k_1}^*\cap\dots \cap    Q_{m,k_m}^* \neq \emptyset$, we can
pick a  cube $R_{k_1, \dots , k_m}$ such that
\begin{equation}\lab{kk3}
Q_{1,k_1}^*\cap \dots \cap Q_{m,k_m}^* \subset
R_{k_1, \dots , k_m} \subset R_{k_1, \dots , k_m}^* \subset
Q_{1,k_1}^{**}\cap \dots \cap Q_{m,k_m}^{**}
\end{equation}
and $|R_{k_1, \dots , k_m}| \ge c |Q_{1,k_1}|$.

Without loss of generality
assume that $Q_{1,k_1}$ has the smallest size among all these cubes.
Since $T$ maps $L^2\times \li\times \dots \times \li$ into $L^2$
we obtain that
\begin{align}\begin{split}\lab{bvnv}
&\int_{R_{k_1,\dots , k_m} } |T(a_{1,k_1}, \dots ,
a_{m,k_m})(x)| \, dx  \\
\le &\bigg(\int_{\rn  } |T(a_{1,k_1}, \dots ,
a_{m,k_m})(x)|^2 \, dx\bigg)^{\f12} |R_{k_1,\dots , k_m}|^{\f12} \\
\le & C(A+B)
|Q_{1,k_1}^{**}|^{\f12}
|Q_{1,k_1}|^{\f12-\f{1}{p_1}} \prod_{j=2}^m |Q_{j,k_j}|^{-\f{1}{p_j}},
\end{split}\end{align}
since $\|a_{1,k_1}\|_{L^2} \le |Q_{1,k_1}|^{\f12-\f{1}{p_1}}$ and
$\|a_{j,k_j}\|_{\li} \le |Q_{j,k_j}|^{ -\f{1}{p_j}}$.
It follows from  (\ref{bvnv})  that
$$
\int_{R_{k_1,\dots , k_m}} |T(a_{1,k_1}, \dots ,
a_{m,k_m}) (x)|  \, dx \le
C(A+B)
|Q_{1,k_1}| \prod_{j=1}^m |Q_{j,k_j}|^{-\f{1}{p_j}}
$$
which combined with $|R_{k_1, \dots , k_m}| \ge c |Q_{1,k_1}|$ gives
\begin{equation}\lab{kk}
\f{1}{|R_{k_1,\dots , k_m}|}\int_{R_{k_1,\dots , k_m}} |T(a_{1,k_1}, \dots ,
a_{m,k_m}) (x)|  \, dx \le
C(A+B)  \prod_{j=1}^m |Q_{j,k_j}|^{-\f{1}{p_j}} .
\end{equation}
We now have the easy estimate
$$
G_1(x) \le \sum_{k_1}\dots \sum_{k_m} |\la_{1, k_1}| \dots |\la_{m, k_m}|
| T(a_{1,k_1}, \dots , a_{m,k_m}) (x)| \chi_{R_{k_1, \dots , k_m}}(x),
$$
and using Lemma \ref{L1},  estimate  (\ref{kk}),
and the last inclusion in (\ref{kk3}) we obtain
\begin{align*}
\big\|G_1\big\|_{\lp} \le & C (A+B) \bigg\|
\sum_{k_1}\dots \sum_{k_m} |\la_{j, k_j}| \dots |\la_{m, k_m}|
\prod_{j=1}^m |Q_{j,k_j}|^{-\f{1}{p_j}} \chi_{Q_{1,k_1}^{**}}
\dots \chi_{Q_{m,k_m}^{**}}\bigg\|_{\lp} \\
\le & C (A+B) \bigg\| \prod_{j=1}^m \bigg(
\sum_{k_j}  |\la_{j, k_j}|
 |Q_{j,k_j}|^{-\f{1}{p_j}} \chi_{Q_{j,k_j}^{**}} \bigg)
 \bigg\|_{\lp} \\
\le & C (A+B) \prod_{j=1}^m  \big \| \big(
\sum_{k_j}  |\la_{j, k_j}|
 |Q_{j,k_j}|^{-\f{1}{p_j}} \chi_{Q_{j,k_j}^{**}} \big)
 \big\|_{L^{p_j}}
\le   C' (A+B) \prod_{j=1}^m  \big \|  f_j \big\|_{H^{p_j}}.
\end{align*}
This proves   (\ref{one}) which combined
with (\ref{two}) completes the proof of the theorem.
\end{proof}

\section{The proof of Lemma \ref{L1}}\label{section3}

It remains to prove Lemma \ref{L1}.  This lemma will be a
consequence of the lemma below. Let $\cd$ be the collection
of all dyadic cubes on $\rn$ and $\cd_j$ be the set of all
cubes in $\cd$ with side length  $l(Q)= 2^{-j}$.

\begin{lemma}\lab{L2}
Suppose $0<p\le 1$.  Then there is a constant $C(p)$
such that  for all finite subsets
 $\mathcal J$ of $\cd$ and all collections    $\{f_Q:\,\, Q\in
\mathcal J\}$   of non-negative integrable functions on $\rn$
with $\text{supp } f_Q \subset Q$ we have
$$
\big\|\sum_{Q\in\mathcal{J}}f_Q\big\|_{\lp} \le \,
C(p)\,  \big\|\sum_{Q\in\mathcal {J }}a_Q\chi_Q \big\|_{\lp},
$$
where
$$
a_Q=|Q|^{-1}\int_Q f_Q(x)\, dx .
$$
\end{lemma}

\begin{proof}
Let us set $\mathcal J_m = \mathcal J \cap \cd_m$ for all $m\in \mathbf Z$.
Given $Q\in \mathcal J$, we define $s(Q)$ to be the unique $m$ such that
$Q\in \mathcal J_m$.
We also set
$$
F= \sum_{Q\in\mathcal{J}}f_Q,\qq
G=\sum_{Q\in\mathcal {J }}a_Q\chi_Q , \qq
G_m= \sum_{k=-\nf}^{m} \sum_{Q\in \mathcal J_k} a_{Q} \chi_Q  .
$$
We now observe that if $Q\in \mathcal J$ and $m\le s(Q)$, then
$G_m$ is constant on $Q$. Therefore for $j\in \mathbf Z$ the
sets below are well-defined
\begin{align*}
\mathcal R_j &=\{Q\in \mathcal J:\,\, G_{s(Q)} \le 2^j \q\text{on $\, Q  $}\},\\
\mathcal R_j' &=\{Q\in \mathcal J:\,\, G_{s(Q)} > 2^j \q\text{on $\, Q\, $ and }
\, G_{s(Q)-1} \le  2^j \q\text{on $\, Q$}\}.
\end{align*}
 For $Q\in \mathcal R_j'$ and any     $t\in Q$, we let
$$
\la_Q= \f{2^j-G_{s(Q)-1}(t)}{G_{s(Q)}(t)-G_{s(Q)-1}(t)}.
$$
Note that $\la_Q$  is a constant since
both functions $G_{s(Q)}$ and $G_{s(Q)-1}$ are constant on $Q$.

We claim that for all $x\in \rn$ we have the identity
\begin{equation}\lab{cc}
 \sum_{Q\in \mathcal R_j}    a_Q \chi_Q(x)   +
\sum_{Q\in \mathcal R_j'} \la_Q a_Q \chi_Q(x)  = \min \big(2^j , G(x) \big)  .
\end{equation}
To prove (\ref{cc}) observe that if $G(x)\le 2^j$, then $\mathcal R_j'=
\emptyset$ and the conclusion easily follows. Otherwise, there is
a smallest $m=m(x)$ such that
$$
2^j <\sum_{k=-\nf}^m \sum_{Q \in \mathcal J_k} a_Q\chi_Q(x).
$$
Then all the cubes that contain $x$ from the collection
$\cup_{k\le m-1}\mathcal  J_k$ belong   to $\mathcal R_j$
and the cube that contains $x$ from $\mathcal  J_m$ belongs to
$\mathcal R_j'$. It follows that
\begin{align*}
 &\sum_{Q\in \mathcal R_j}    a_Q \chi_Q(x)   +
\sum_{Q\in \mathcal R_j'} \la_Q a_Q \chi_Q(x)  \\
=& G_{m-1}(x) + \sum_{Q\in \mathcal R_j'} (2^j-G_{m-1}(x)) \chi_Q(x) =2^j,
\end{align*}
since the last sum has only one term. This proves (\ref{cc}).
Next we set
$$
F_j = \sum_{Q\in \mathcal R_j} f_Q + \sum_{Q\in \mathcal R_j'} \la_Q f_Q .
$$
Then using (\ref{cc}) we obtain
\begin{align}\begin{split}\lab{453}
\int_{\rn}  F_j(x)\, dx  =  &\int_{\rn}  \bigg(
 \sum_{Q\in \mathcal R_j}
  a_Q \chi_Q(x)   +  \sum_{Q\in \mathcal R_j'} \la_Q a_Q \chi_Q(x)   \bigg)
dx \\ = &\int_{\rn}  \min \big(2^j , G(x) \big)\, dx   .
\end{split}\end{align}
Now $F_j-F_{j-1}$ is supported on $\{G>2^{j-1}\}$. H\"older's inequality
and (\ref{453}) give
\begin{align*}
\int_{\rn} (F_j(x)-F_{j-1}(x))^p\,  dx &\le |\{G>2^{j-1}\}|^{1-p}
\bigg(\int_{\rn} \!\min \big(2^j , G(x) \big)\, dx \bigg)^{p}\\
&\le 2^{jp} |\{G>2^{j-1}\}|.
\end{align*}
Summing the above over all $j$
and using the fact that
$$
F(x) = \sum_{j\in \mathbf Z} \big(F_{j}(x)-F_{j-1}(x)\big)
$$
we obtain the required estimate
$$
\int_{\rn} (F (x) )^p\,  dx \le
\sum_{j\in \mathbf Z} 2^{jp} |\{G>2^{j-1}\}| \le C(p)^p
\int_{\rn} (G (x) )^p\,  dx ,
$$
where the last inequality follows by summation
by parts.
\end{proof}

Having established Lemma \ref{L2}, we now proceed to the proof of
Lemma \ref{L1}.

\begin{proof}
Given the  cubes $\{Q_k\}_{k=1}^K$ we can find
a finite collection of dyadic cubes $\{Q_{kj}\}_{j=1}^{m_k}$
with
$$
  l(Q_k)  \le l(Q_{kj}) \le 2l(Q_k)
$$
and
\begin{equation}\lab{ok}
Q_k \subset \bigcup\limits_{j=1}^{m_k} Q_{kj} \subset Q_k^*,
\end{equation}
where $m_k\le 2^n$.
We apply Lemma \ref{L2} to the functions
$\{g_k \chi_{Q_{kj}}\}_{1\le k \le K}^{1\le j \le m_k}$.
(We collapse terms when the same dyadic cube is used twice). We obtain
\begin{equation}\lab{last}
\big\| \sum_{k=1}^K g_k \big\|_{\lp} \le
\big\|
\sum_{k=1}^K  \sum_{j=1}^{m_k} g_{k} \chi_{Q_{kj}}
\big\|_{\lp} \le  C (p) \big\|
\sum_{k=1}^K  \sum_{j=1}^{m_k} b_{kj} \chi_{Q_{kj}}
\big\|_{\lp} ,
\end{equation}
where $b_{kj} = |Q_{kj}|^{-1} \dint_{Q_{kj}} g_k(x)\, dx \le
|Q_{k}|^{-1} \dint_{Q_{k}} g_k(x)\, dx$. Inserting this estimate in
(\ref{last}) gives
$$
\big\| \sum_{k=1}^K g_k \big\|_{\lp} \le  C(p)
\bigg\| \sum_{k=1}^K \bigg(|Q_{k}|^{-1} \dint_{Q_{k}} g_k(x)\, dx
\bigg)   \sum_{j=1}^{m_k}   \chi_{Q_{kj}}  \bigg\|_{\lp},
$$
and the required conclusion follows from the last inclusion in (\ref{ok}).
\end{proof}

\section{Related results and comments}\label{section4}
We note that Theorem \ref{main} can be extended to the
case when some $p_j$'s
are bigger than $1$ and the remaining $p_j$'s  are less than or equal to
$1$. We have the following:

\begin{theorem}\lab{main2}
Let   $1<q_1,\dots , q_m,q<\nf$ be fixed indices satisfying   (\ref{ind1})
and let $0<p_1,\dots , p_m,p<\nf$ be any real numbers satisfying (\ref{ind2}).
Suppose that $K$ satisfies (\ref{assump}) for all $|\al |\le N$
where $N$ is sufficiently large.
Let $T$ be related to $K$ as in (\ref{defT}) and assume that $T$
admits an extension that
maps $L^{q_1}(\rn)\times \dots \times L^{q_m}(\rn)$ into
$L^q(\rn)$ with norm $B$.  Then
$T$ extends to a bounded operator from
$H^{p_1}(\rn)\times \dots \times H^{p_m}(\rn)$ into
$L^p(\rn)$ (we set $H^{p_j}=L^{p_j}$ when ${p_j}>1$), which satisfies
the norm estimate
$\|T\|_{H^{p_1} \times \dots \times H^{p_m}\to L^p} \le C (A+B)  $
for some   constant $C=C(n,p_j,q_j)$. ($A$ is as in (\ref{defa}).)
\end{theorem}

\begin{proof}
We discuss the multilinear interpolation needed to prove
this theorem for all indices $0<p_j<\nf$.
Theorem \ref{main2} is valid when all
the $p_j$'s satisfy $1<p_j<\nf$ as proved in \cite{GT}.
In Theorem \ref{main} we considered the case when all $0<p_j\le 1$.

We now fix indices $0<p_j<\nf$ so that some
of them are bigger than $1$ and some of them are less than or equal to
$1$. We pick $\ve>0$ and $\la>0$ so that
$$
0<\ve <\min \big(\tfrac{1}{m}, \tfrac{1}{p_1} , \dots , \tfrac{1}{p_m} \big),
\qqq
\la > \big(\min \big(  \tfrac{1}{p_1} , \dots , \tfrac{1}{p_m} \big)
-\ve\big)^{-1}.
$$
Using that $T$ is bounded from $L^{q_1}\times  \dots \times L^{q_m}$
into $L^q$ with norm at most $B$, it follows from
 \cite{GT}  that
\begin{equation}\lab{int1}
T:\,\, L^{1/\ve}\times \dots \times L^{1/\ve}\to L^{1/m\ve}
\end{equation}
with norm at most a constant multiple of $A+B$.
Now define $s_j$ be setting
\begin{equation}\lab{aol}
\tfrac{1}{s_j} = \la \big( \tfrac{1}{p_j}-\ve\big) +\tfrac{1}{p_j}.
\end{equation}
Then it is easy to see that $0<s_j<1$ for all $1\le j\le m$ and by
Theorem \ref{main}  we have
\begin{equation}\lab{int2}
T:\,\, H^{s_1}\times \dots \times H^{s_m}\to L^{s},
\end{equation}
with norm at most a constant multiple of $A+B$,
where $1/s=1/s_1+\dots + 1/s_m$.
Here we need (\ref{assump}) with $N=[n(1/s-1)]$.
Identity (\ref{aol}) gives
$$
\tfrac{1}{p_j}= \tfrac{\theta}{s_j}+\tfrac{1-\theta}{1/\ve}
$$
where $\theta =  (\la+1)^{-1}$. Interpolating between
(\ref{int1}) and (\ref{int2}) we obtain
that
$$
T:\,\, [L^{1/\ve},H^{s_1}]_\theta \times \dots \times
[L^{1/\ve},H^{s_m}]_\theta \to [L^{1/m\ve},L^{s }]_\theta = L^p,
$$
where $1/p=1/p_1+\dots + 1/p_m$. But
$[L^{1/\ve},H^{s_j}]_\theta = L^{p_j}$ if $p_j>1$ or $H^{p_j}$ if
$p_j\le 1$  and the required conclusion follows (see e.g. \cite{JJ}).
\end{proof}

We note that mapping into the Hardy space $H^p$ instead of $L^p$
is hopeless even in the translation invariant case unless some
further cancellation is imposed. We refer to \cite{GT2} and
\cite{cdm} for results of this sort.  Both references deal with
bilinear operators but the  techniques can be adapted to give
similar results  for   $m$-linear operators as well.

We now discuss some analogous results for the maximal
singular integral operator defined by
$$
T_*(f_1,\dots , f_m)(x)= \sup_{\de>0} |T_\de (f_1,\dots , f_m)(x)|,
$$
where $T_\de$ are the smooth truncations of $T$ given by
$$
T_\de (f_1,\dots , f_m)(x)   =
\int_{\rn}
K_\de(x,y_1,\dots , y_m) f_1(y_1)\dots f_m(y_m) \, dy_1\dots dy_m.
$$
Here $K_\de(x,y_1,\dots , y_m) = \eta \big(
\sqrt{|  x \!-\! y_1|^2 +\dots + |  x \!-\!  y_m|^2} /\de   \big)
K (x,y_1,\dots , y_m)$ and $\eta$ is a smooth
function on $\rn$ which vanishes  in a neighborhood of the
origin and is equal to $1$ outside a larger neighborhood of the origin.

It is proved in \cite{GT3} that the sublinear operator $T_*$ satisfies
similar boundedness
estimates as $T$. We have the following result regarding $T_*$.

\begin{theorem}\lab{main3}
Under the same hypotheses are Theorem \ref{main2},
$T_*$ maps the product
$H^{p_1}(\rn)\times \dots \times H^{p_m}(\rn)$ boundedly into
$L^p(\rn)$ , and satisfies
the norm estimate
$\|T_*\|_{H^{p_1} \times \dots \times H^{p_m}\to L^p} \le C (A+B)  $
for some   constant $C=C(n,p_j,q_j)$.  As usually, we set
$H^{p_j}=L^{p_j}$ when ${p_j}>1$.
\end{theorem}

\begin{proof}
The proof is similar to that for $T$.
First we consider the case where all the $p_j$'s are less than or equal
to one.  It follows from    \cite{GT3} that $T_*$
is bounded
on the same range as $T$ with bound at most a multiple of $A+B$.
Thus the estimates in case 1 follow as before.
Next observe that the kernels $K_\de$ satisfy  (\ref{assump})
uniformly in $\de>0$. Hence
the estimates in case 2 for $K $
equally apply to $K_\de $   uniformly in $\de>0$ and the
same conclusion follows.

The remainder of the argument is then similar.  One treats the
multilinear maps
$$
T_{\delta_1,\ldots,\delta_N}(f_1,\dots , f_m)(x)=  \{T_{\de_k}
(f_1,\dots
, f_m)(x)\}_{k=1}^N,
$$
as maps $T_{\delta_1,\ldots,\delta_N}: H^{s_1}\times \dots \times
H^{s_m}\to L^{s}(\ell_{\infty}^N),$ for any finite set
$\delta_1,\ldots,\delta_N>0$ and uses complex interpolation as before.
\end{proof}

The first author would like to thank Xuan Thinh Duong for
his hospitality in Sydney where  part of this work was conceived.

\end{document}